\documentclass{article}

\usepackage{textcomp}
\usepackage{geometry}
\usepackage{xcolor}
\usepackage{float}
\usepackage{mathrsfs}
\usepackage{amsmath}
\usepackage{amssymb}
\usepackage{graphicx}
\usepackage{color}
\usepackage{cite}
\makeatletter

\floatstyle{ruled}
\newfloat{algorithm}{tbp}{loa}
\providecommand{\algorithmname}{Algorithm}
\floatname{algorithm}{\protect\algorithmname}

\newtheorem{thm}{\protect\theoremname}
\newtheorem{problem}[thm]{\protect\problemname}
\newtheorem{rem}[thm]{\protect\remarkname}

\newcommand{\prox}{{\rm prox}}

\providecommand{\problemname}{Problem}
\providecommand{\remarkname}{Remark}
\providecommand{\theoremname}{Theorem}

\makeatother

\providecommand{\problemname}{Problem}
\providecommand{\remarkname}{Remark}
\providecommand{\theoremname}{Theorem}

\begin{document}
\title{Optimal Mass Transport of Nonlinear Systems under Input and Density
Constraints}

\author{Dongjun Wu, Anders Rantzer 
 \thanks{*The authors are supported by the Excellence Center ELLIIT and the European Research Council (Advanced Grant 834142).}
	\thanks{The authors are with Department of Automatic Control, Lund
	University, Box 118, SE-221 00 Lund, Sweden {\tt\small
	dongjun.wu, anders.rantzer@control.lth.se}.}
	}
	\date{}

\maketitle
\begin{abstract}
We investigate optimal mass transport problem of affine-nonlinear dynamical systems
with input and density constraints. Three algorithms are proposed 
to tackle this problem, including
two Uzawa-type methods and a splitting algorithm based on the Douglas-Rachford algorithm. 
Some preliminary simulation results are presented to demonstrate the effectiveness of our approaches.
\end{abstract}

\section{Introduction}
In recent years, there has been a growing interest in leveraging optimal mass transport techniques in various
control applications, including multi-robotic swarm \cite{krishnan2018distributed},
mean-field
control \cite{fornasier2014mean,fornasier2014meansparse}, and optimal steering
\cite{chen2016optimal,chen2015optimal}. These applications typically aim to control an initial density,
representing the probability distribution of multi-agent systems, towards a target density while minimizing
associated costs.

The problem is not new. 
In the 2000s, important theoretical advancements were made utilizing optimal control theory and partial
differential equation techniques. Bernard and Buffoni, for instance, characterized the existence and
regularity of optimal plans using Hamilton-Jacobi equation \cite{bernard2007optimal}. Agrachev and Lee
further explored optimal transport of non-holonomic systems employing Pontryagin maximum principle
\cite{agrachev2009optimal}, see also an earlier related work by Ambrosio et al. \cite{ambrosio2004optimal}.
More recently, Ghoussoub et al. extended these findings to the free terminal time case
\cite{ghoussoub2018optimal}.

Another significant line of research stems from the seminal work of Benamou and Brenier, who introduced a
computational fluid dynamics approach to solving the optimal transport problem
\cite{benamou2000computational}. Their method has been particularly successful in numerical computations due
to the convex nature of the problem formulation, enabling the development of efficient algorithms, see e.g.,
\cite{papadakis2014proximal}. It is also worth mentioning that \cite{rantzer2001dual} applied similar techniques to convexify the search for
stabilizing nonlinear controllers. Building upon this foundation, Chen et al. systematically
investigated optimal transport and steering for linear dynamical systems, under which explicit solutions
are sometimes available \cite{chen2016optimal,chen2015optimal}.
More recently, Elamvazhuthi et al. explored optimal
transport of affine nonlinear systems \cite{elamvazhuthi2023dynamical}. 
Nevertheless,
numerical algorithms capable of handling input and density constraints
were not proposed. Similarly, Kerrache and Nakauchi
extended Benamou and Brenier's work to include density constraints \cite{kerrache2022constrained}, but did not
consider dynamical systems nor input constraints.

In this paper, we address the challenge of optimal transport of affine nonlinear systems under both input and
density constraints. The contribution is three algorithms which solve the problem successfully.
The first two are based on the recent work \cite{kerrache2022constrained} and the third one is based
on the celebrated Douglas-Rachford algorithm.

\section{Problem formulation}

We study optimal mass transport (OMT) over the following nonlinear
dynamical system:
\begin{equation}
\dot{x}=f(x,t)+B(x,t)u\label{sys:nlti}
\end{equation}
where $x\in X\subseteq\mathbb{R}^{n}$ denotes the state,
$u=(u_{1},\cdots,u_{m})\in U\subseteq\mathbb{R}^{r}$ represents the input
and $B(x,t)={\rm col}(b_{1}(x,t),\cdots,b_{r}(x,t))$. The vector fields $f$
and $b_{i}$s are assumed to be continuously differentiable and $B(x,t)$
has full rank for all $x$ and $t$.
The set $X$ and $U$ are both compact. 
We call $u(\cdot)$ an admissible control as long as $u(t)\in U$
for all $t\ge0$.
The problem of optimal transport amounts
to finding a control policy $u$ to steer an initial density of agents with dynamics
\eqref{sys:nlti} to a target density with minimum cost.




The main interest of this paper is to incorporate constraints
in the optimal mass transport
of nonlinear dynamical systems \eqref{sys:nlti}. We consider
two types of constraints. The first type is constraints on the density
flow $\rho(t,\cdot)$ of the agents. For example, in scenarios where agents
must avoid spedific obstacles at time $t$, $\rho(t,x)$ should vanish
at those spatial points. Similarly, if the density is not allowed 
to exceed a certain level,
it should be bounded accordingly. We use the function
$I(\rho)$ to represent the density constraint. For example when $\rho(t,x)\le h(x)$
for all $t$ and $x$ for a given function $h$, then $I(\rho)=1_{\{\phi:\;\phi\le h\}}(\rho)$
where 
\[
1_{A}(x)=\begin{cases}
0, & x\in A\\
+\infty, & \text{otherwise}
\end{cases}
\]

Meanwhile, we impose constraints on the input, which is implicitly
encoded in the admissible control region $U$. This consideration
is relevant when dealing with physical applications where
the control actions are subject to constraints. 
While a recent paper \cite{elamvazhuthi2023dynamical}
considered the case $u\in U$ where $U$ is convex, 
the algorithm proposed therein appears to be applicable
only in constraint-free scenarios, i.e., when $U=\mathbb{R}^{m}$.

Motivated by these considerations, we propose to investigate the following constrained
OMT problem:
\begin{problem}
Given an initial density $\rho_{0}\in\mathscr{P}(X)$
of homogeneous agents, each governed by dynamics (\ref{sys:nlti}), determine
an admissible control law $u$ such that the density of the agents
at the terminal time $T$ is $\rho_{T}\in\mathscr{P}(X)$, while minimizing the
cost functional 
\begin{equation}
J(\rho,u)=\frac{1}{2}\int_{0}^{T}\int_{X}\rho(t,x)|u(t,x)|^{2}dxdt+I(\rho)\label{cost:J(u)} 
\tag{J}
\end{equation}
under the PDE constraints 
\begin{equation}
\left\{ \begin{aligned} & \partial_{t}\rho+\nabla\cdot(\rho(f+Bu))=0\\
 & \rho(0,\cdot)=\rho_{0}\\
 & \rho(T,\cdot)=\rho_{T}
\end{aligned}
\right.\label{eq:PDE-density}
\end{equation}
\end{problem}

\section{Main Results}

We make the following standing assumptions:

\textbf{Assumption 1: }The system \eqref{sys:nlti} is controllable.

\textbf{Assumption 2: } Both the control region $U$ and the density 
constraint function $I(\rho)$ are convex.

\subsection{Uzawa Type Algorithm: A Direct Method}
Following \cite{benamou2000computational} and the recent work \cite{elamvazhuthi2023dynamical},
we introduce a Lagrangian multiplier $\phi\in C^{1}(\mathbb{R}\times X;\mathbb{R})$
to handle the PDE constraints (\ref{eq:PDE-density}). Assume that
$\rho$ is compactly supported in the interior of
$X$, then applying integration
by parts yields
\begin{align*}
 & \inf_{\rho,u}J(\rho,u)\\
= & \inf_{\rho,m}\sup_{\phi}\int\frac{|m|^{2}}{2\rho}-\rho\partial_{t}\phi-Bm\cdot\nabla\phi dxdt-G(\phi)+I(\rho)
\end{align*}
where we have denoted $m=\rho u$ and
\begin{equation}
G(\phi)=\int_{X}\phi(0,x)\rho_{0}(x)-\phi(T,x)\rho_{T}(x)dx.\label{eq:G(phi)}
\end{equation}

Following \cite{benamou2000computational}, we transform the 
term $|m|^{2}/2\rho$ into a more tractable form: 
\[
\frac{|m|^{2}}{2\rho}=\sup_{a,b}a\rho+m\cdot b,\quad(a,b)\in K
\]
where 
\begin{equation}
K:=\left\{ (a,b):a+\frac{|b|^{2}}{2}\le0\right\} \label{K}
\end{equation}
This transformation has two merits: 1) it ensures
$\rho$ to be non-negative and 2) the cost $a\rho+m\cdot b$
becomes linear in the optimization variable $(\rho,m)$.

Denoting $q=(a,b)$, and $\mu=(\rho,m)$, the original problem 
is reformulated as follows:
\begin{align*}
 & \inf_{\rho,u}J(\rho,u)\\
= & \inf_{\mu}\sup_{q,\phi} -I_{K}(q)-G(\phi)+I(\rho)+\left\langle \mu,q-P\nabla_{t,x}\phi\right\rangle _{H}
\end{align*}
where 
\[
P(x)=\begin{bmatrix}1 & f(x)^{\top}\\
0 & B(x)^{\top}
\end{bmatrix}\in\mathbb{R}^{(r+1)\times(n+1)}.
\]
and $\left\langle \cdot,\cdot\right\rangle _{H}$ denotes the inner
product in the Hilbert space 
\begin{equation}
H=:L^{2}([0,T]\times X;\mathbb{R}^{r+1}),\label{eq:spaceV}
\end{equation}
that is,
\[
\left\langle \mu_{1},\mu_{2}\right\rangle _{H}=\int_{[0,T]\times X}\mu_{1}\cdot\mu_{2}dxdt.
\]

Remember that we still have a constraint on $u$, whereas our optimization
variables are $\mu,q,\phi$. Thus it remains to rephrase the input constraint
in the language of the new variables. To that end, recalling that
$U$ is convex, we can express it as $U=\cap_{\alpha}\{z\in\mathbb{R}^{r}:a_{\alpha}^{\top}z+b_{\alpha}\le0\}$
for suitable vectors $(a_{\alpha},b_{\alpha})_{\alpha\in\mathcal{A}}$.
But this is equivalent to \footnote{The formula should be understand in a point-wise sense at non-vanishing
points of $\rho$.} 
\begin{equation*}
(\rho,m)\in\mathcal{U}=\bigcap_{\alpha}\{(\rho,m)\in\mathbb{R}_{+}\times\mathbb{R}^{r}:a_{\alpha}^{\top}m+b_{\alpha}\rho\le0\}.\label{eq:rho-m-cvx}
\end{equation*}
and $\mathcal{U}$ is again a convex set. Using $I(\mu)$ to include both
the input constraints and the density constraint
$I(\rho)$, we now propose the final formulation of the problem under
consideration: 
\begin{align}
 & \sup_{\mu}\inf_{q,\phi}L(\mu,q,\phi)  \label{supinf-P}\tag{P1}\\
= & \sup_{\mu}\inf_{q,\phi} I_{K}(q)+G(\phi)-I(\mu)+\left\langle \mu,P\nabla_{t,x}\phi-q\right\rangle _{H}.  \nonumber
\end{align}

For this problem to be well-defined, we need a Hilbert space
in which the function $\phi$ lives.
Since the cost involves time and space derivatives,
the natural candidate is the Sobolev space $H^{1}([0,T]\times X)$.
However, for technical reasons,
we shall, as
in \cite{benamou2000computational,elamvazhuthi2023dynamical}, impose
an extra condition that $\int\phi(t,x)dxdt=0$. Thus we define 
\[
V:=\left\{ \phi\in H^{1}([0,T]\times X);\;\int\phi dxdt=0\right\} 
\]
with inner product 
\[
\left\langle \phi_{1},\phi_{2}\right\rangle _{V}=\int\phi_{1}\phi_{2}+\nabla_{t,x}\phi_{1}\cdot\nabla_{t,x}\phi_{2}dxdt.
\]

This type of problem has been previously studied in \cite{kerrache2022constrained}, 
where it was tackled using Uzawa-type algorithms, 
see Algorithm \ref{alg:cons-opt-1}. However, the algorithm
therein requires the operator $E=P\nabla $ to be injective and
has closed range, which however, is not satisfied in our case.
To remedy this, 
let $B_{2}(x,t)\in\mathbb{R}^{n\times(n-r)}$ be such that $\tilde{B}(x,t):={\rm col}(B,B_{2})$
is invertible. 
Introduce an additional control input $v$ so that the system
(\ref{sys:nlti}) becomes 
\[
\dot{x}=f(x)+\tilde{B}(x,t)\begin{bmatrix}u\\
v
\end{bmatrix}.
\]
while imposing the linear constraint $v=0$.
Now we can replace $P$ in \eqref{supinf-P}
by the invertible matrix $\tilde{P}(x,t)\in\mathbb{R}^{(n+1)\times(n+1)}$:
\[
\tilde{P}(x,t)=\begin{bmatrix}1 & f^{\top}\\
0 & \tilde{B}^{\top}
\end{bmatrix}(x,t)
\]
and that the operator $E=\tilde{P}\nabla$ becomes injective
and has closed range.
The
algorithm is summarized in Algorithm \ref{alg:cons-opt-1} where $E\phi=P\nabla\phi$.

\begin{algorithm}[ht]
\caption{Direct approach}
\label{alg:cons-opt-1}

\textbf{Initialization}: $(p^{0},b^{0},\mu^{1},\nu^{0},\eta^{0})$
given arbitrarily.

\textbf{Repeat:} 
\begin{enumerate}
\item Compute $\phi^{n}$ by solving 
\begin{align*}
\phi^{n}=\arg\min_{\phi}G(\phi) & +\left\langle \nu^{n-1},E\phi-p^{n-1}\right\rangle _{H}\\
 & +\frac{r}{2}\|E\phi-p^{n-1}\|_{H}^{2}
\end{align*}
\item $p^{n}=p^{n-1}-\rho_{r}(\mu^{n}-\nu^{n-1}+r(p^{n-1}-B\phi^{n})$; 
\item $\nu^{n}=\nu^{n-1}+\rho_{s}(E\phi^{n}-p^{n-1}-s(\nu^{n-1}-\mu^{n}))$; 
\item Compute $q^{n}$ by solving 
\[
q^{n}=\arg\min_{q\in K}\left|b^{n-1}+\frac{\eta^{n-1}}{r}-q\right|^{2}
\]
pointwisely. 
\item $b^{n}=b^{n-1}-\rho_{r}(\eta^{n-1}-\mu^{n}+r(b^{n-1}-q^{n})$; 
\item $\eta^{n}=\eta^{n-1}+\rho_{s}(b^{n-1}-q^{n}-s(\eta^{n-1}-\mu^{n}))$; 
\item Compute $\mu^{n+1}$ by solving 
\[
\hspace{-7mm}\mu^{n+1}=\arg\min_{\mu}s\left|\mu-[\nu^{n}+\eta^{n}+(p^{n}-b^{n})/s]/2\right|+I(\mu)
\]
pointwisely. 
\end{enumerate}
\end{algorithm}

In Algorithm \ref{alg:cons-opt-1}, 
several auxiliary optimization variables, namely $p,b,\nu$ and $\eta$
along with some positive optimization parameters, namely $r,s,\rho_{r},\rho_{s}$
need to be introduced. It has been demonstrated in \cite{kerrache2022constrained}
that the algorithm converges, at least weakly, to a saddle point $(\phi^{*},q^{*},\mu^{*})$\footnote{The convergence of the auxiliary variables are 
incidental to us.}
under the following additional assumptions: 
\begin{enumerate}
\item The auxiliary cost function
\begin{align*}
\mathcal{L}_{r,s}( & \phi,q,p,b,\mu,\nu,\eta)=F(q)+G(\phi)-I(\mu)\\
 & +\left\langle \mu,p-b\right\rangle _{H}+\left\langle \nu,P\nabla\phi-p\right\rangle _{H}+\left\langle \eta,b-q\right\rangle _{H}\\
 & +\frac{r}{2}\|P\nabla\phi-p\|_{H}^{2}+\frac{r}{2}\|b-q\|_{H}^{2}\\
 & -\frac{s}{2}\|\mu-\nu\|_{H}^{2}-\frac{s}{2}\|\mu-\eta\|_{H}^{2}
\end{align*}
admits a saddle point $(\phi^{*},q^{*},p^{*},b^{*},\mu^{*},\nu^{*},\eta^{*})$. 
\item The optimization parameters satisfy\footnote{Such set of parameters always exist, see \cite{kerrache2022constrained}.}
\begin{align*}
2s-\rho_{r}-\rho_{s}s^{2}-|\rho_{r}r-\rho_{s}s| & >0,\\
2r-\rho_{r}r^{2}-\rho_{s}-|\rho_{r}r-\rho_{s}s| & >0.
\end{align*}
\end{enumerate}

\begin{rem}
    The existence of saddle points of the optimization problem 
    \eqref{supinf-P} is a 
    delicate issue and can potentially be analyzed using
    the same techniques as in the work \cite{hug2020convergence}.
    
\end{rem}

The most pivotal part of this algorithm lies in step 1), which amounts
to solving the following partial differential equation: 
\begin{equation}
-r\nabla_{t,x}\cdot(\tilde{P}^{\top}\tilde{P}\nabla_{t,x}\phi^{n})
=\nabla_{t,x}\cdot(\tilde{P}^{\top}\nu^{n-1}-r\tilde{P}^{\top}p^{n-1})\label{eq:PDE-step1}
\end{equation}
with pure Neumann boundary conditions 
\begin{align*}
r(\tilde{P}^{\top}\tilde{P}\nabla_{t,x}\phi^{n})_{0}(0,\cdot) & =\rho_{0}-(\tilde{P}^{\top}\nu^{n-1})_{0}(0,\cdot)\\
 & \quad\quad+r(\tilde{P}^{\top}p^{n-1})_{0}(0,\cdot)\\
r(\tilde{P}^{\top}\tilde{P}\nabla_{t,x}\phi^{n})_{0}(T,\cdot) & =\rho_{T}-(\tilde{P}^{\top}\nu^{n-1})_{0}(T,\cdot)\\
 & \quad\quad+r(\tilde{P}^{\top}p^{n-1})_{0}(T,\cdot).
\end{align*}
The presence of the matrix function $\tilde{P}$ in (\ref{eq:PDE-step1}) 
can lead to potential challenges for the PDE: 1) it may
possibly degrade the well-posedness of the equation and may 2) have worse 
regularity properties
compared to the PDEs encountered in \cite{benamou2000computational}
and \cite{elamvazhuthi2023dynamical}. In response, we propose
an alternative approach which only involves solving a standard Poisson
equation.

\subsection{Uzawa Type Algorithm: An indirect method}

We introduce a new vector field $v=f(x,t)+B(x,t)u$. 
Then the cost function (\ref{cost:J(u)})
can be recast as:
\[
J=\frac{1}{2}\int\rho|B^{\dag}(v-f)|^{2}dxdt+I(\rho)
\]where $B^\dag$ is the Moore-Penrose inverse of $B$.

Denote $m=\rho v$ as the new momentum variable, we can
rewrite $J$ as 
\[
J=\int\frac{|B^{\dag}m|^{2}}{2\rho}+\frac{\rho}{2}|B^{\dag}f|^{2}-B^{\dag}m
\cdot B^{\dag}fdxdt+I(\rho)
\]
After introducing the set of admissible momentum variables: 
\begin{align*}
\mathfrak{M}_{{\rm ad}}=\{m=  \rho(f+Bu)|\
  u(t)\in U, \; \forall t\in[0,T]\}.
\end{align*}
Problem 1 can be reformulated as the optimization problem 
\[
\min_{\rho\ge0,\;m\in\mathfrak{M}_{{\rm ad}}}J(\rho,m).
\]

Let $\phi\in C^{1}(\mathbb{R}\times X;\mathbb{R})$ be a Lagrangian multiplier.
The cost (\ref{cost:J(u)}) can be recast as 
\begin{align*}
J(\rho,& m,\phi)  = 
\int_{0}^{T}\int_{X}\frac{|B^{\dag}m|^{2}}{2\rho}+\frac{\rho}{2}|B^{\dag}f|^{2}\\
 & \quad-B^{\dag}m\cdot B^{\dag}f+\phi\left\{ \frac{\partial\rho}{\partial t}+{\rm div}\left\{ v\right\} \right\} dxdt + I(\rho)  \\
 & =\int_{0}^{T}\int_{X}\frac{|B^{\dag}m|^{2}}{2\rho}+\frac{\rho}{2}|B^{\dag}f|^{2}-B^{\dag}m\cdot B^{\dag}f \\
 & \quad-\rho\frac{\partial\phi}{\partial t}-\nabla\phi\cdot mdxdt-G(\phi)
 + I(\rho)
\end{align*}
where $G(\phi)$ is as (\ref{eq:G(phi)}). Remembering the fact that 
\[
\frac{|B^{\dag}m|^{2}}{2\rho}=\sup_{a,b\in K}(a\rho+b\cdot B^{\dag}m),
\]
with $K$ identical to \eqref{K}, Problem 1 can further be transformed to
\begin{align*}
\inf_{\rho,\;m\in\mathfrak{M}_{{\rm ad}}}\sup_{q,\phi}\int\rho(\alpha-\partial_{t}\phi) & +m\cdot(\rho-\nabla_{x}\phi)dxdt\\
 & -G(\phi)+I(\rho)+I_{K_{f}}(q)
\end{align*}
in which 
\[
K_{f}=:\{(\alpha,\beta):\alpha+ B^\dag \beta \cdot B^\top f+\frac{1}{2}
|B^{\top}\beta|^{2}\le0\}.
\]
Indeed, 
\begin{align*}
 & \inf_{\rho\ge0,\;m\in\mathfrak{M}_{{\rm ad}}}\sup_{\phi}J(\rho,m,\phi)\\
= & \inf_{\rho\ge0,m\in\mathfrak{M}_{{\rm ad}}}\sup_{(a,b)\in K,\;\phi}\int\{\rho(a+\frac{1}{2}|B^{\dag}f|^{2}-\partial_{t}\phi)\\
 & +\left.m\cdot((B^{\dag})^\top b- 
 (B^{\dag})^\top B^{\dag}f-\nabla_{x}\phi)\right\} dxdt\\
 & -G(\phi)+I(\rho)
\end{align*}
Letting 
\[
\alpha=a+\frac{1}{2}|B^{\dag}f|^{2},\quad
\beta=(B^{\dag})^\top b-(B^{\dag})^\top B^{\dag}f
\]
yields the above formulation.

It still remains to remove the constraints
$m\in\mathfrak{M}_{{\rm ad}}$. Denote $\mu=(\rho,m)$ and let $I(\mu)$
represent the overall input constraints
and density constraints $I(\rho)$ . The problem is thus reformulated as:
\begin{equation}
-\sup_{\mu}\inf_{\phi,q}I_{K_{f}}(q)+G(\phi)-I(\mu)+\left\langle \mu,\nabla_{t,x}\phi-q\right\rangle _{H}\label{-maxmin} \tag{P2}
\end{equation}

The optimization problem (\ref{-maxmin}) again mirrors the
structure in \cite{kerrache2022constrained}. The key difference
lies in $K_{f}$, which is now dependent on the drift vector $f$ and
the input matrix $B(x,t)$.
Thus $K_f$ encapsulates the complete dynamical 
system information.

Consequently, we can apply Algorithm \ref{alg:cons-opt-1} again where
$E$ simply becomes $\nabla_{t,x}$ -- an injective operator with closed range -- 
ensuring the convergence of the algorithm under the aforementioned conditions.
The algorithm is outlined in Algorithm \ref{alg:cons-opt-2}.

\begin{algorithm}[ht]
\caption{Indirect approach}
\label{alg:cons-opt-2}

\textbf{Initialization}: $(p^{0},b^{0},\mu^{1},\nu^{0},\eta^{0})$
given arbitrarily.

\textbf{Repeat:} 
\begin{enumerate}
\item Compute $\phi^{n}$ by solving the PDE 
\[
-r\Delta_{t,x}\phi^{n}=\nu^{n-1}-rp^{n-1}
\]
with pure Neumann boundary conditions 
\begin{align*}
r\partial_{t}\phi^{n}(0,\cdot) & =\rho_{0}-\nu_{0}^{n-1}(0,\cdot)+rp_{0}^{n-1}(0,\cdot)\\
r\partial_{t}\phi^{n}(T,\cdot) & =\rho_{T}-\nu_{0}^{n-1}(T,\cdot)+rp_{0}^{n-1}(T,\cdot).
\end{align*}
\item $p^{n}=p^{n-1}-\rho_{r}(\mu^{n}-\nu^{n-1}+r(p^{n-1}-\nabla_{t,x}\phi^{n})$; 
\item $\nu^{n}=\nu^{n-1}+\rho_{s}(\nabla_{t,x}\phi^{n}-p^{n-1}-s(\nu^{n-1}-\mu^{n}))$; 
\item Compute $q^{n}$ by solving 
\[
q^{n}=\arg\min_{q\in K_{f}}\left|b^{n-1}+\frac{\eta^{n-1}}{r}-q\right|^{2}
\]
point-wisely where 
\[
K_{f}=:\{(\alpha,\beta):\alpha+ B^\dag \beta \cdot B^\top f+\frac{1}{2}
|B^{\top}\beta|^{2}\le0\}.
\]
\item $b^{n}=b^{n-1}-\rho_{r}(\eta^{n-1}-\mu^{n}+r(b^{n-1}-q^{n})$; 
\item $\eta^{n}=\eta^{n-1}+\rho_{s}(b^{n-1}-q^{n}-s(\eta^{n-1}-\mu^{n}))$; 
\item Compute $\mu^{n+1}$ by solving 
\[
\hspace{-7mm}\mu^{n+1}=\arg\min_{\mu}s\left|\mu-[\nu^{n}+\eta^{n}+(p^{n}-b^{n})/s]/2\right|+I(\mu)
\]
pointwisely. 
\end{enumerate}
\end{algorithm}

\begin{rem}
It is noteworthy that although we intially required
the inverse of the matrix $B$ -- which could be unbounded or discontinuous --
to reframe the optimization problem, in the end, the need vanishes entirely.
\end{rem}

\subsection{Douglas-Rachford algorithm}

Douglas-Rachford algorithm (DR) is a splitting algorithm tailored for
addressing the optimization problem \cite{douglas1956numerical}:
\begin{equation}
\inf_{x}F_{1}(x)+F_{2}(x)\label{eq:DR-form} \tag{DR}
\end{equation}
This algorithm is useful when the proximal operators of
the convex functions $F_{1}$ and $F_{2}$
are readily computable. Given an initial pair $(x^{0},\bar{x}^{0}),$ the DR algorithm
iterates through the following steps:
\[
\begin{cases}
\hat{x}^{k}=2x^{k}-\bar{x}^{k}\\
\bar{x}^{k+1}=\overline{x}^{k}+\prox_{F_{1}}(\hat{x}^{k})-x^{k}\\
x^{k+1}=\prox_{F_{2}}(\bar{x}^{k+1}).
\end{cases}
\]

To reframe Problem 1 into the form of \eqref{eq:DR-form}, we work with the
original formulation rather than adding Lagragian multipliers as in the previous
two algorithms.
For clarity, we make the following assumption.

\textbf{Assumption 3: } The state space is $X=[0,1]^{n}$.
The density constraint is $g(t,x)\le\rho(t,x)\le h(t,x)$
for some non-negative functions $g$ and $h$, while the input constraints
take the form $|u_{i}|\le k_{i}$ for some non-negative constants
$k_{i}$ for all $i$.

Under this assumption, the constraints can be expressed as 
\[
R\mu\le\gamma
\]
for some constant matrix $R$ and a fixed vector-valued function $\gamma$.
Now introduce a linear (unbounded) operator $A$ 
\[
A\mu=\begin{bmatrix}\partial_{t}\rho+\nabla\cdot m\\
\rho(0,\cdot)\\
\rho(T,\cdot)
\end{bmatrix},
\]
and a vector $\theta=[0,\rho_{0},\rho_{T}]^{\top}$. 
Then Problem 1 can be stated as
\begin{equation}
\inf_{\rho,m}M(\mu)+I_{\theta}(A\mu)+I_{C_{\gamma}}(R\mu)\label{eq:DR_eq0}
\tag{P3}
\end{equation}
where 
\[
M(\mu)=\int\frac{|B^{\dag}m|^{2}}{2\rho}+\frac{\rho}{2}|B^{\dag}f|^{2}-B^{\dag}m\cdot B^{\dag}fdxdt
\]
and $C_{\gamma}=\{\xi:\xi\le\gamma\}$.

The cost function of \eqref{eq:DR_eq0}
contains three terms. 
Various splitting algorithms exists to tackle optimization problems
of this type, namely $\min_{x} f(x)+g(x)+h(Ax)$,
see e.g., \cite{condat2013primal, vu2013splitting}. Typically, these algorithms require
$f$ to be differentiable and $A$ to be bounded, which do not fit into 
our setting. One may potentially modify the Hilbert space, say e.g., Sobelev space
$H^1$ for $\mu$. But this leads to additional complexities. 
Thus we employ DR algorithm which converges under very mild condition: 
when \eqref{eq:DR_eq0} admists a minimizer
and that $F$ and $G$ are convex.

To transform it into the standard form \eqref{eq:DR-form},
we introduce two additional variables $\eta$ and $\xi$ 
so that \eqref{eq:DR_eq0} becomes:
\[
\inf_{\mu,\eta,\xi}\left\{ M(\mu)+I_{\theta}(\eta)+I_{C_{\gamma}}(\xi)\right\} +I_{\{A\mu=\eta,R\mu=\xi\}}(\mu,\eta,\xi)
\]
Denoting
\begin{align*}
F(\mu,\eta,\xi) & =M(\mu)+I_{\theta}(\eta)+I_{C_{\gamma}}(\xi)\\
N(\mu,\eta,\xi) & =I_{\{A\mu=\eta,R\mu=\xi\}}(\mu,\eta,\xi)
\end{align*}
reduces the problem to standard form \eqref{eq:DR-form}.

The function $F$ has the form $\sum_{i}f_{i}(x_{i})$, whose proximal
can be separated: 
\begin{equation}
\prox_{F}(\mu,\eta,\xi)=\prox_{M}(\mu)+\prox_{I_{\theta}}(\eta)+\prox_{I_{C_{\gamma}}}(\xi)
\end{equation}
where each proximal operator can be readily computed. 
The proximal of $N$
is a projection to a linear subspace which is 
the most involved part of the algorithm. We proceed to calculate the
proximal operators one by one.

1) $\prox_{M(\mu)}$. This amounts to solving the point-wise optimization
problem: 
\begin{align*}
\prox_{M}(\mu) & =\arg\min_{\tilde{\mu}}\frac{1}{2}|\tilde{\mu}-\mu|^{2}\\
 & \quad+\frac{|B^{\dag}\tilde{m}|^{2}}{2\tilde{\rho}}+\frac{\tilde{\rho}}{2}|B^{\dag}f|^{2}-B^{\dag}\tilde{m}\cdot B^{\dag}f.
\end{align*}

2) $\prox_{I_{\theta}}$. Since $\theta$ is a singleton, $\prox_{I_{\theta}}(\eta)=\theta$
for all $\eta$. 

3) $\prox_{I_{C_{\gamma}}}$. This is a point-wise optimization on the space
$[0,T]\times X$: $
\prox_{I_{C_{\gamma}}}(\xi)=\arg\min_{\tilde{\xi}\le\gamma}\frac{1}{2}|\tilde{\xi}-\xi|^{2}
$.

4) $\prox_{N}$. We show in the appendix that this amounts to solving
a partial differential equation. 

The complete algorithm is outlined in Algorithm \ref{alg:cons-opt-2}.


\subsection{Numerical Example: Double integrator}

\begin{figure}[ht]
\begin{centering}
\includegraphics[scale=0.2]{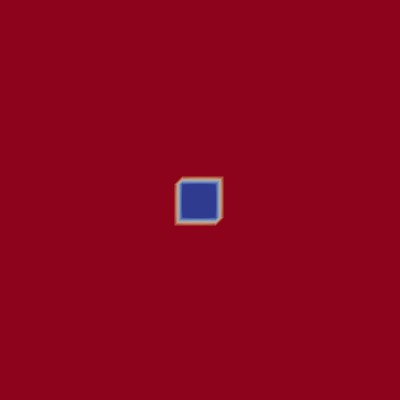}
\par\end{centering}
\caption{Density constraint.}\label{fig:density_const}

\end{figure}

Consider the double integrator $ \dot{x}_{1} =x_{2}$, $\dot{x}_{2}  =u$,
on $X=[0,1]\times[0,1]$ with initial density 
$\rho_{0}(x)$ which is equal to $10$ on the region
$(x_{1}-0.25)^{2}+(x_{2}-0.4)^{2}<0.15$ and $0$ otherwise.
The target density $\rho_{T}(x)$ is equal to $10$
on the region $(x_{1}-0.75)^{2}+(x_{2}-0.6)^{2}<0.15$ and $0$ otherwise.
Assume $T=1$. We consider a density constraints $\rho\le h$ which is equal to $10$
on the region $ \max\{|x_{1}-0.5|,|x_{2}-0.5|\}\le0.05$ and $0$ otherwise:
this represents an obstacle in the center of the space, see Fig. \ref{fig:density_const}.

\begin{figure}[ht]
\begin{centering}
\includegraphics[scale=0.2]{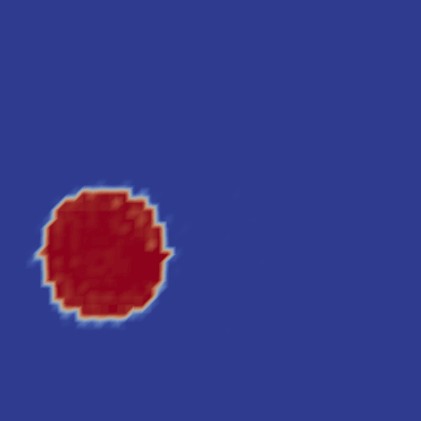}
\includegraphics[scale=0.2]{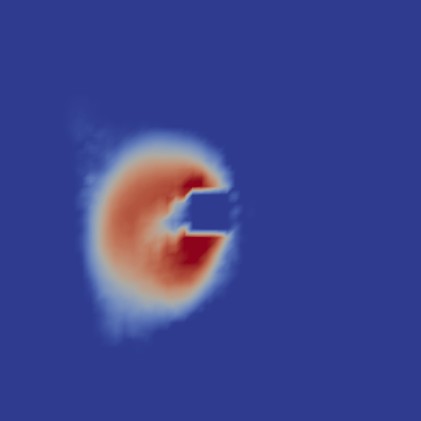}
\includegraphics[scale=0.2]{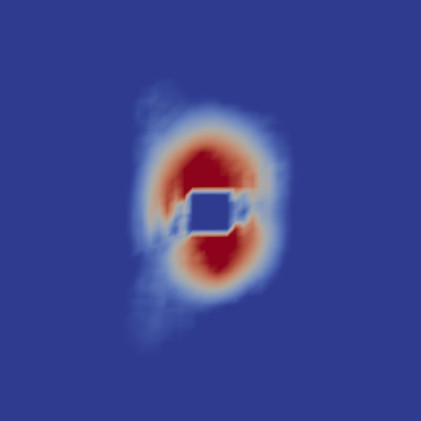}
\par\end{centering}
\vspace{1.3mm}
\begin{centering}
\includegraphics[scale=0.2]{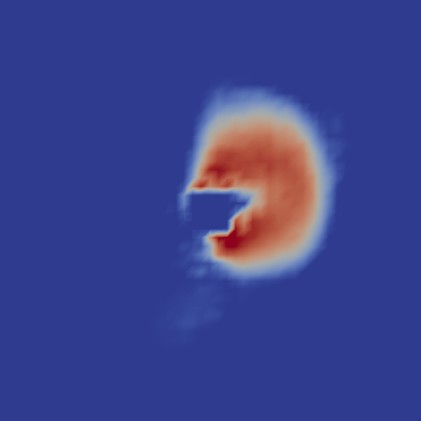}
\includegraphics[scale=0.2]{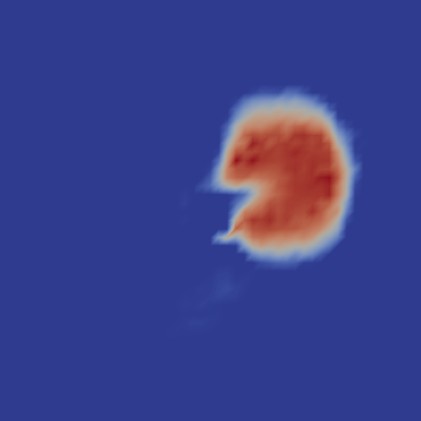}
\includegraphics[scale=0.2]{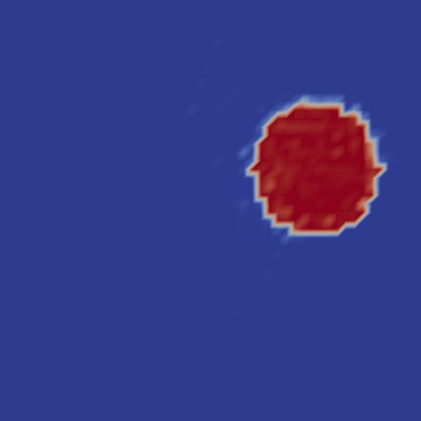}
\par\end{centering}
\caption{Density flow under constraints $|u|\le 20$}
\label{fig:rho_u20}
\end{figure}

\begin{figure}[ht]
\begin{centering}
\includegraphics[scale=0.2]{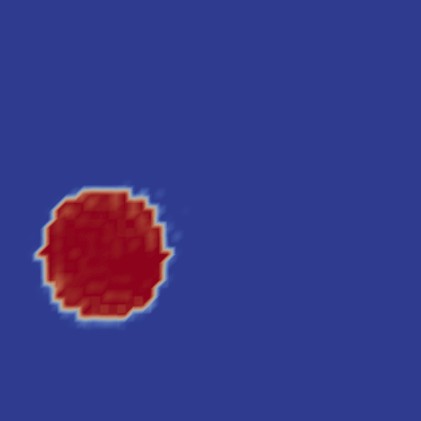}
\includegraphics[scale=0.2]{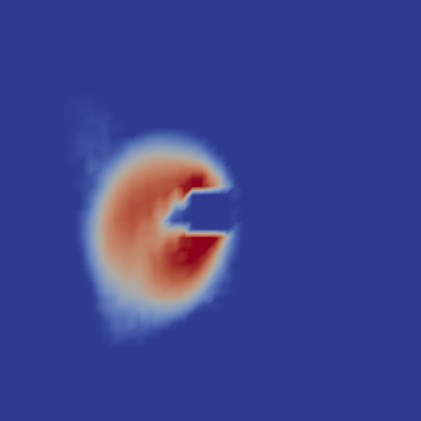}
\includegraphics[scale=0.2]{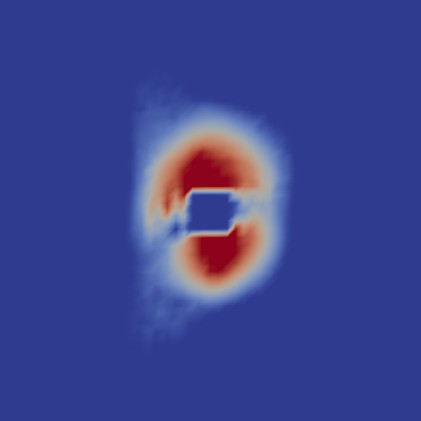}
\par\end{centering}
\vspace{1.3mm}
\begin{centering}
\includegraphics[scale=0.2]{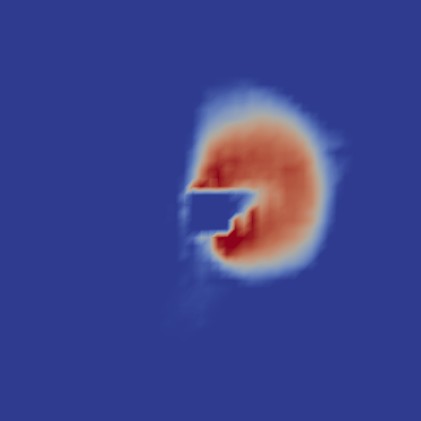}
\includegraphics[scale=0.2]{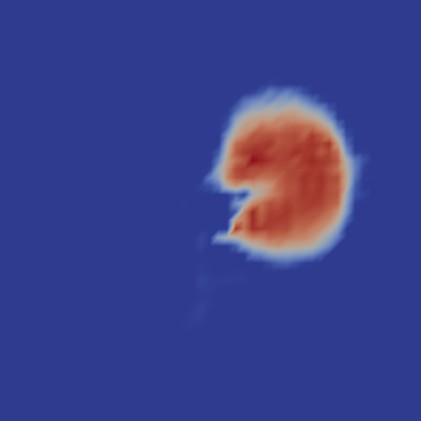}
\includegraphics[scale=0.2]{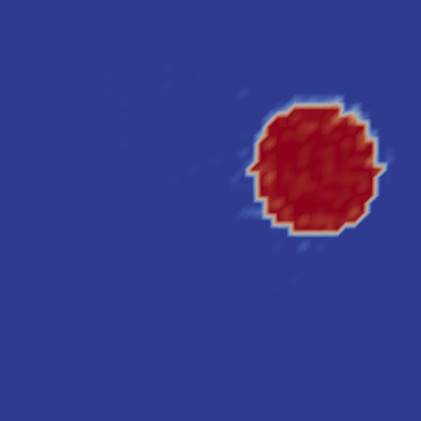}
\par\end{centering}
\caption{Density flow under constraints $|u|\le 1$.}
\label{fig:rho_u1}
\end{figure}

The optimal transport plans under constraints $|u|\le 20$ and
$|u|\le 1$ are shown in Fig. \ref{fig:rho_u20} and \ref{fig:rho_u1} respectively.
The simulation is done using Algorithm \ref{alg:cons-opt-2}.
We can see that the density flows are comparable 
under the two constraints.

\begin{figure}

\begin{centering}
\includegraphics[scale=0.2]{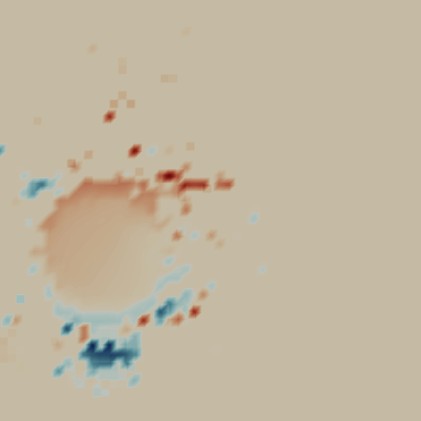}
\includegraphics[scale=0.2]{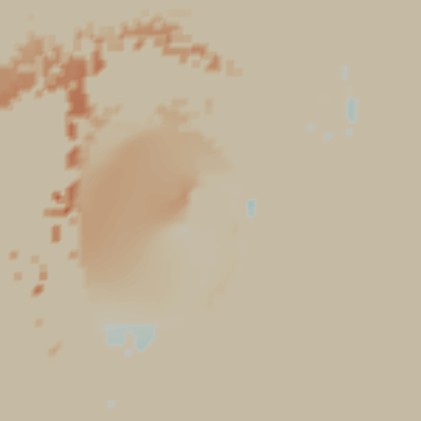}
\includegraphics[scale=0.2]{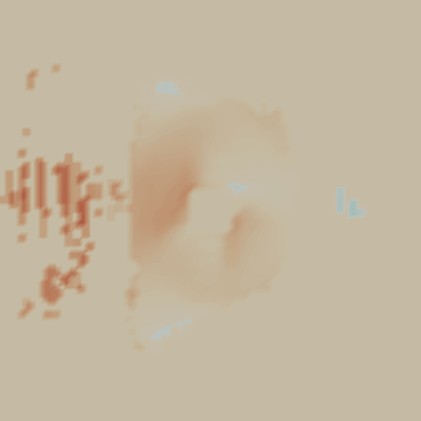}
\par\end{centering}
\vspace{1.3mm}
\begin{centering}
\includegraphics[scale=0.2]{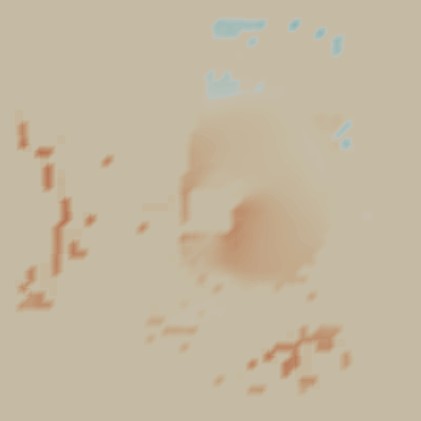}
\includegraphics[scale=0.2]{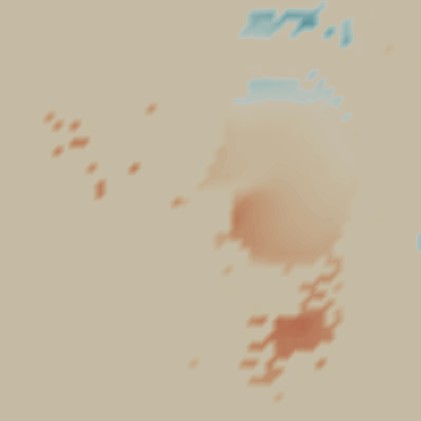}
\includegraphics[scale=0.2]{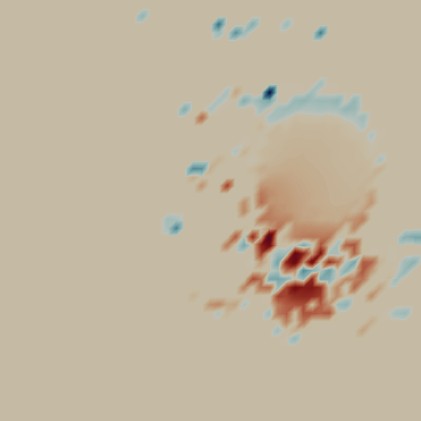}
\par\end{centering}

\caption{Input on space $|u|\le 20$.}\label{fig:input_u20}

\end{figure}

\begin{figure}[ht]
\begin{centering}
\includegraphics[scale=0.2]{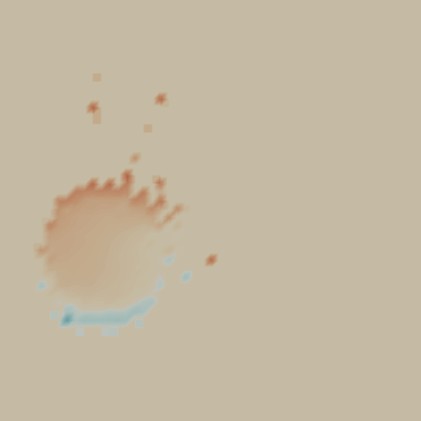}
\includegraphics[scale=0.2]{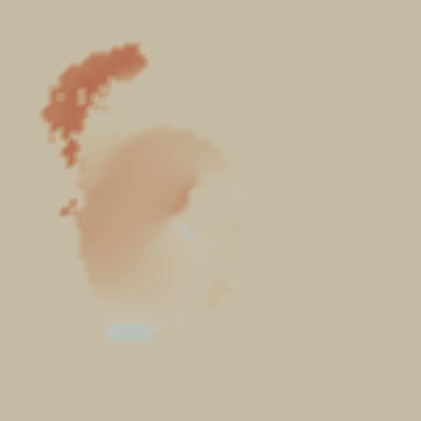}
\includegraphics[scale=0.2]{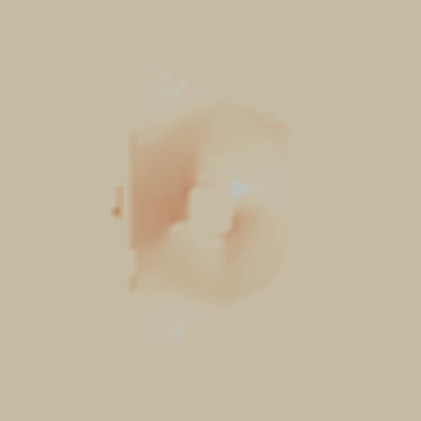}
\par\end{centering}
\vspace{1.3mm}
\begin{centering}
\includegraphics[scale=0.2]{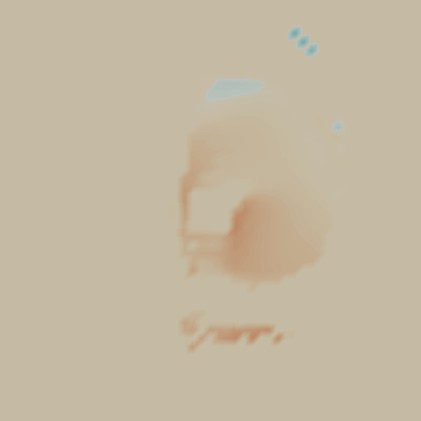}
\includegraphics[scale=0.2]{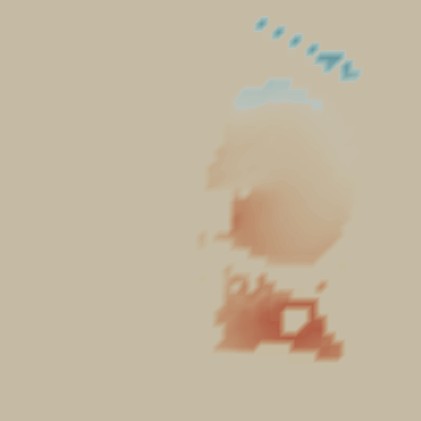}
\includegraphics[scale=0.2]{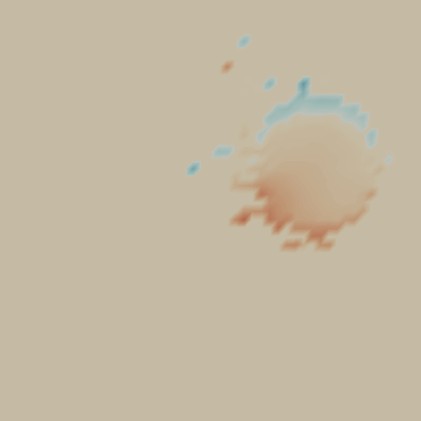}
\par\end{centering}
\caption{Input on space, $|u|\le 1$. }\label{fig:input_u1}
\end{figure}

Fig \ref{fig:input_u20} and \ref{fig:input_u1} show the corresponding control input 
in space and time -- the warmer the color, the higher the control actions.
We can observe that, the control constraints somehow smooth the distribution of the control
input in space and that when there is no input constraints, 
large control actions may appear which is not allowed in practice.

\section{Conclusion}
In this paper, we have studied optimal mass transport of affine nonlinear
control systems under input and density constraints. The proposed algorithms
has guaranteed convergence under mild conditions. Future direction includes
extending these results to systems on large scale graphs.

\section{Acknowledgement}
We would like to thank Professor Pontus Giselsson for insightful discussions and valuable
feedback during the preparation of this manuscript.

We have used dolfinx \cite{barrata2023dolfinx} to solve the PDEs in this paper.
And the computations were enabled by resources provided by the Swedish National Infrastructure for Computing (SNIC) at the PDC Center for High Performance Computing, KTH Royal Institute of Technology, partially funded by the Swedish Research Council through grant agreement no. 2018-05973.

\section{Appendix}

\subsection*{Computation of $\prox_{N}(\mu,\eta,\xi)$}

By definition of the proximal operator
\begin{align}
& {\rm prox}_{N}(\mu,\eta,  \xi)  \nonumber \\
= & \inf_{\tilde{\mu},\tilde{\eta},\tilde{\xi}}\int\frac{1}{2}|\tilde{\mu}
-\mu|^{2}+\frac{1}{2}|\tilde{\xi}-\xi|^{2} 
 +\frac{1}{2}|\tilde{\eta}_{0}-\eta_{0}|^{2}dxdt \nonumber \\
& \quad +\int\frac{1}{2}|\tilde{\eta}_{1}-\eta_{1}|^{2}+\frac{1}{2}|\tilde{\eta}_{2}-\eta_{2}|^{2}dx\label{appen:proxJ}
\end{align}
under the constraints 
\[
\begin{cases}
\nabla_{t,x}\cdot\tilde{\mu}=\tilde{\eta}_{0}\\
\tilde{\rho}(0,\cdot)=\tilde{\eta}_{1},\tilde{\rho}(T,\cdot)=\tilde{\eta}_{2}\\
R\tilde{\mu}=\tilde{\xi}
\end{cases}
\]
We substitute the last two constraints into \eqref{appen:proxJ},
and then we are left with two optimization variables, namely, $\tilde{\mu}$
and $\tilde{\eta}_{0}$. Introduce a Lagrangian multiplier $\phi(t,x)$
and define $\phi_{Q}(t,x)=\phi\left(Q^{-1}\begin{bmatrix}t\\
x
\end{bmatrix}\right)$ for the first constraint, where
\[
Q = (I+R^\top R)^{-1}.
\]

The problem then becomes
\begin{align*}
\inf_{\tilde{\mu},\tilde{\eta}_{0}}\sup_{\phi} & \int\frac{1}{2}|\tilde{\mu}-\mu|^{2}+\frac{1}{2}|R\tilde{\mu}-\xi|^{2}+\frac{1}{2}|\tilde{\eta}_{0}-\eta_{0}|^{2}dxdt\\
 & +\int\frac{1}{2}|\tilde{\rho}(0,\cdot)-\eta_{1}|^{2}+\frac{1}{2}|\tilde{\rho}(T,\cdot)-\eta_{2}|^{2}dx\\
 & -\int Q^{-1}\nabla\phi\cdot\tilde{\mu}+\phi\tilde{\eta}_{0}dxdt\\
 & +\int\phi_{Q}(T,\cdot)\tilde{\rho}(T,\cdot)-\phi_{Q}(0,\cdot)\tilde{\rho}(0,\cdot)dx
\end{align*}

\begin{algorithm}[ht]
\caption{Douglas-Rachford based Algorithm\protect}
\label{alg:cons-opt-3}

\textbf{Initialization}: $(\mu^{0},\eta^{0},\xi^{0},\overline{\mu}^{0},\overline{\eta}^{0},\overline{\xi}^{0})$
given arbitrarily.

\textbf{Repeat:}
\begin{enumerate}
\item Calculate $\hat{\mu}^{k},\hat{\eta}^{k},\hat{\xi}^{k}$:
\begin{align*}
\hat{\mu}^{k} & =2\mu^{k}-\overline{\mu}^{k},\\
\hat{\eta}^{k} & =2\eta^{k}-\overline{\eta}^{k},\\
\hat{\xi}^{k} & =2\xi^{k}-\overline{\xi}^{k}
\end{align*}
\item Calculate $\mu^{k+1},\eta^{k+1},\xi^{k+1}$:
\begin{enumerate}
\item $\mu^{k+1}=\mu^{k}+\Phi(\hat{\mu}^{k})-\bar{\mu}^{k}$, where $\Phi(\hat{\mu}^{k})$
is the solution to the finite dimensional nonlinear equation (only
the positive root of $\tilde{\rho}$ is relevant):
\begin{equation}
\left\{ \begin{aligned} & \rho-\hat{\rho}^{k}-\frac{|B^{\dag}m|^{2}}{2(\hat{\rho}^{k})^{2}}+\frac{1}{2}|B^{\dag}f|^{2}=0\\
 & m-\hat{m}^{k}+\frac{(B^{\dag})^\top B^{\dag}m}{\rho}-
 (B^{\dag})^\top B^{\dag}f=0
\end{aligned}
\right.\label{eq:NL_eq}
\end{equation}
with dependent variables $(\rho,m)$.
\item $\eta^{k+1}=\eta^{k}+\theta-\bar{\eta}^{k}$.
\item Solve $\xi^{k+1}=\arg\min_{\xi\le\gamma}\frac{1}{2}|\xi-\hat{\xi}^{k}|^{2}$
pointwisely.
\end{enumerate}
\item Find $\phi^{k}$ by solving the PDE:
\[
-\Delta_{t,x}\phi^{k}+\phi^{k}=\nabla\cdot(\mu^{k+1}+R^{\top}\xi^{k+1})-\eta_{0}^{k+1}
\]
with boundary conditions:
\begin{equation}
\left\{
\begin{aligned}
&\left.\partial_{t}\phi^{k}\right|_{t=0}-\phi^{k}(0,\cdot)
=\eta_{1}^{k+1} \\ 
& \hspace{3cm} -(\mu^{k+1}+R^{\top}\xi^{k+1})_{0}(0,\cdot)\\
& \left.\partial_{t}\phi^{k}\right|_{t=T}+\phi^{k}(T,\cdot)
=\eta_{2}^{k+1} \\
& \hspace{3cm}-(\mu^{k+1}+R^{\top}\xi^{k+1})_{0}(T,\cdot)\\
& \left.\partial_{i}\phi^{k}\right|_{x_{i}=0,1}=-(R^{\top}\xi^{k+1})_{i}|_{x_{i}=0,1},\quad\forall i
\end{aligned}\label{PDE-DR:bd}
\right.
\end{equation}
where $\partial_{i}\phi^{k}=\partial\phi^{k}/\partial x_{i}$. 
\item Update $\overline{\mu}^{k+1},\overline{\eta}^{k+1},\overline{\xi}^{k+1}$:
\begin{align*}
\bar{\mu}^{k+1} & =Q^{-1}(\nabla\phi^{k}+\mu^{k+1}+R^{\top}\xi^{k+1})\\
\bar{\eta}_{0}^{k+1} & =\phi^{k}+\eta_{0}^{k+1}\\
\overline{\eta}_{1}^{k+1} & =\rho^{k+1}(0,\cdot)-\phi^{k}(0,\cdot)\\
\bar{\eta}_{2}^{k+1} & =\rho^{k+1}(T,\cdot)+\phi^{k}(T,\cdot)\\
\bar{\xi}^{k+1} & =R\bar{\mu}^{k+1}
\end{align*}
\end{enumerate}
\end{algorithm}

{\noindent}from which the KKT condition can be derived:
\begin{align*}
\tilde{\mu}-\mu+R^{\top}(R\tilde{\mu}-\xi)-Q^{-1}\nabla\phi & =0\\
\tilde{\eta}_{0}-\eta_{0}-\phi & =0\\
\tilde{\rho}(0,\cdot)-\eta_{1}-\phi_{Q}(0,\cdot) & =0\\
\tilde{\rho}(T,\cdot)-\eta_{2}+\phi_{Q}(T,\cdot) & =0
\end{align*}
Solving $\phi$ from the KKT condition and plugging into the cost,
we get
\begin{align*}
\sup_{\phi} & \int\frac{1}{2}|b+\nabla\phi-\mu|^{2}+\frac{1}{2}|R(b+\nabla\phi)-\xi|^{2}+\frac{1}{2}\phi^{2}dxdt\\
 & +\int\frac{1}{2}\phi_{Q}(0,\cdot)^{2}+\frac{1}{2}\phi_{Q}(T,\cdot)^{2}dx\\
 & +\int-Q^{-1}\nabla\phi\cdot(b+\nabla\phi)-\phi(\eta_{0}+\phi)dxdt\\
 & +\int\phi_{Q}(T,\cdot)(\eta_{2}-\phi_{Q}(T,\cdot))-\phi_{Q}(0,\cdot)(\eta_{1}+\phi_{Q}(0,\cdot))dx
\end{align*}
where $b=Q(\mu+R^{\top}\xi)$. To find the optimal $\phi$, it is
sufficient to equate the variation of the cost at $\phi$ to $0$.
That is, for all $\varphi$, there must hold
\begin{align*}
0=- & \int\nabla\varphi\cdot(\mu+R^{\top}\xi+Q^{-1}\nabla\phi)-\varphi(\phi+\eta_{0})dxdt\\
+\int & \varphi(T,\cdot)(\eta_{2}-\phi_{Q}(T,\cdot))-\varphi(0,\cdot)(\eta_{1}+\phi_{Q}(0,\cdot))dx
\end{align*}
Applying integration by parts yields
\begin{align*}
 & \int\varphi\{(\nabla\cdot(\mu+R^{\top}\xi+\nabla\phi_{Q})-(\phi+\eta_{0})\}dxdt\\
+ & \int\varphi(0,\cdot)[(\mu+R^{\top}\xi)_{0}(0,\cdot)+\partial_{t}\phi_{Q}(0,\cdot)]dx\\
+ & \int-\varphi(T,)[(\mu+R^{\top}\xi)_{0}(T,\cdot)+\partial_{t}\phi_{Q}(T,\cdot)]dx\\
+ & \sum_{i=1}^{n}\varphi|_{x_{i}=0}[(\mu+R^{\top}\xi)_{i}+\partial_{i}\phi_{Q})]_{x_{i}=0}\\
+ & \sum_{i=1}^{n}\varphi|_{x_{i}=1}[(\mu+R^{\top}\xi)_{i}+\partial_{i}\phi_{Q})]_{x_{i}=1}\\
+ & \int\varphi(T,\cdot)[\eta_{2}-\phi_{Q}(T,\cdot)]-\varphi(0,\cdot)[\eta_{1}+\phi_{Q}(0,\cdot)]dx
\end{align*}
Since $\varphi$ is arbitrary, we immediately get
\begin{align*}
\Delta\phi_{Q}+\nabla\cdot(\mu+R^{\top}\xi)-(\phi_{Q}+\eta_{0})=0\\
\partial_{t}\phi_{Q}(0,\cdot)+(\mu+R^{\top}\xi)_{0}(0,\cdot)-\eta_{1}(\cdot)-\phi_{Q}(0,\cdot)=0\\
-\partial_{t}\phi_{Q}(T,\cdot)-(\mu+R^{\top}\xi)_{0}(T,\cdot)+\eta_{2}(\cdot)-\phi_{Q}(T,\cdot)=0\\
\partial_{i}\phi_{Q}|_{x_{i}=0,1}+(R^{\top}\xi)_{i}|_{x_{i}=0,1}=0
\end{align*}
Once $\phi_{Q}$ has been found, we can solve for the optimal solution
by employing the KKT condition.

\bibliographystyle{plain}
\bibliography{OT,IEEEabrv}

\end{document}